\documentstyle[12pt,amscd]{amsart}
\newtheorem{Theorem}{Theorem}[section]
\newtheorem{Lemma}[Theorem]{Lemma}
\newtheorem{Corollary}[Theorem]{Corollary}
\newtheorem{Proposition}[Theorem]{Proposition}
\newtheorem{Remark}[Theorem]{Remark}
\newtheorem{Example}[Theorem]{Example}

\def\depth{\operatorname{depth}}
\def\reg{\operatorname{reg}}
\def\Ext{\operatorname{Ext}}

\def\Deg{\operatorname{Deg}}

\def\sk{\smallskip\par}
\def\mm{{\mathfrak m}}

\def\pp{{\mathfrak p}}

\def\FilM{{\mathbb M}}
\def\FilN{{\mathbb N}}

\begin{document}
\title{Dependence of Hilbert coefficients}
\thanks{This is a corrected version of the original paper  published in {\small manuscripta math. {\bf 149},  235 - 249 (2016)}.  In this version we include a Corrigendum and give a small modification of the proof of Theorem 2.4.\\  \   \\
Both authors were partially supported by NAFOSTED (Project  101.01-2011.48). The paper was completed during the stay of the second author at the Vietnam Institute for Advanced Study in Mathematics.\\ {\it 2000 Mathematics Subject Classification:}  Primary 13D40, Secondary 13A30.\\ {\it Key words and phrases:}  Castelnuovo-Mumford regularity, associated graded module,   good filtration, Hilbert coefficient.}

\maketitle
\begin{center}
LE XUAN DUNG\\
Department of  Natural Science,  Hong Duc University\\
307 Le Lai, Thanh Hoa, Vietnam\\
E-mail: lxdung27@@gmail.com\\ [10pt]
and\\[10pt]
LE TUAN HOA \\
 Institute of Mathematics Hanoi (VAST)\\  18 Hoang Quoc Viet 10307 Hanoi, Vietnam \\
 E-mail: lthoa@@math.ac.vn

  \end{center}

\begin{abstract} Let $M$ be a finitely generated module of dimension $d$  and depth $t$ over a Noetherian local ring ($A, \mm$) and $I$ an $\mm$-primary ideal. In the main result it is shown that the last $t$ Hilbert coefficients $e_{d-t+1}(I,M),..., e_d(I,M)$ are bounded below and above in terms of the first $d-t+1$ Hilbert coefficients $e_0(I,M),...,e_{d-t}(I,M)$ and $d$.
\end{abstract}

\date{}
\section*{Introduction} \sk
Let $M$ be a finitely generated module of dimension $d$ over a Noetherian local ring ($A, \mm$) and $I$ an $\mm$-primary ideal. The Hilbert-Samuel function $H_{I,M} (n) = \ell(M/I^{n+1}M)$ agrees with the Hilbert-Samuel polynomial $P_{I,M}(n)$ for $n \gg 0$ and we may write
$$P_{I,M} (n) = e_0(I,M){n+d \choose d} - e_1(I,M){n+d-1 \choose d-1} + \cdots + (-1)^d e_d(I,M).$$
The numbers $e_0(I,M), e_1(I,M), ...,e_d(I,M)$ are called the Hilbert coefficients of $M$ with respect to $I$. 

The Hilbert-Samuel function and the Hilbert-Samuel polynomial give a lot of information on $M$. Therefore, it is of interest to study properties of Hilbert coefficients. Assume that $A$ is a Cohen-Macaulay ring and $M$ is a Cohen-Macaulay $A$-module. Then Northcott \cite{No} and Narita \cite{Na} showed that $e_1(I,A) \geq 0$ and $e_2(I,A) \geq 0$, respectively. Note that already  $e_3(I,A)$ maybe negative. Later, Rhodes \cite{Rh} showed that the above results also hold for  good $I$-filtrations  of submodules of  $M$. Moreover, Kirby and Mehran \cite{KM} were able to show that $e_1(I,M)\leq {e_0(I,M) \choose 2}$ and $e_2(I,M)\leq {e_1(I,M) \choose 2}$. Subsequently these results were improved by several authors. How about  the other coefficients? In 1997, Srinivas and Trivedi \cite{ST2} and Trivedi \cite{Tr1} obtained a surprising result, stating that   all $|e_i(I,A)|, i \geq 1,$  are bounded by a function depending only on $e_0(I,A)$ and $d$. 

What happens for non-Cohen-Macaulay modules? Inspired by the previously  mentioned result of Srinivas and Trivedi and of  Trivedi \cite{Tr2}, Rossi-Trung-Valla \cite{RTV} showed that all $|e_i(I,A)|$,  are bounded by  functions depending on the so-called extended degree $\Deg(I,A)$ and $d$. These results were extended to modules in \cite{L2} and \cite{DH}.  It is also worth to know, that when $I$ is a parameter ideal in a generalized Cohen-Macaulay ring, there is a uniform bound for all  $|e_i(I,A)|,\ i\ge 1$, which does not depend on the choice of $I$, see \cite{GO}. However from all these  results one cannot deduce further relations between Hilbert coefficients.

Using a bound on the Castelnuovo-Mumford regularity in terms of Hilbert coefficients given in \cite[Theorem 2]{Tr1} one can immediately see  that $(-1)^{i-1}e_i(I,A)$ is bounded above by a (complicated and implicit) function depending only on $e_0(I,A),...,e_{i-1}(I,A)$  and $i$, for all $i\ge 1$.  An explicit bound will be given in Theorem \ref{B0}. However,  even in the case $d = 1$ an easy example shows  that $|e_1(I,A)|$ is in general not bounded in terms of $e_0(I,A)$. So, it is natural to ask: how many Hilbert coefficients are enough to  be taken such that they completely bound the absolute values of all other ones?   The main result of this paper is to show that the  first $d-t+1$ Hilbert coefficients have this property, where $t= \depth M$ (see Theorem \ref{B6} and Corollary \ref{B7}). 
As a consequence, we can show that there is only a finite number of Hilbert-Samuel functions if $e_0(I,M), e_1(I,M), ...,e_{d-t}(I,M)$ and $d$ are fixed (see Theorem \ref{B8}).

In fact, we will deal with a more general situation, namely with good $I$-filtrations $\FilM$. In this case our bounds also involve the so-called reduction number $r(\FilM)$. Our approach is somewhat similar to that of \cite{ST2, Tr1} and \cite{RTV}, in the sense that we use  the Castelnuovo-Mumford regularity $\reg(G(\FilM))$ of the associated module $G(\FilM)$ of $\FilM$ to bound the Hilbert coefficients (see Proposition \ref{NB5}).  Then one has to bound $\reg(G(\FilM))$ in terms of  the first $d-t+1$ Hilbert coefficients.   In  order to do that, in Section \ref{CMHil}, using \cite[Theorem 2]{Tr1} we  first give a bound for  $\reg(G(\FilM))$  in terms of all Hilbert coefficients (see Theorem \ref{NB3}).  Then, combining some idea  developed in the proof of   \cite[Theorem 3.3]{RTV},  and  refined in \cite[Theorem 4.4]{L1} and \cite[Theorem 1.5]{DH}, with bounding the length of certain Artinian modules (see Lemma \ref{B4}), we  show in the same section that  already the first $d-t+1$ Hilbert coefficients are enough to bound $\reg(G(\FilM))$ (see Theorem \ref{B5}).  The relations among the Hilbert coefficients  are given in the last section (Theorem \ref{B0} and Theorem \ref{B6}). Finally, we would like  to remark, that bounds established in this  paper are huge functions. Therefore instead of seeking better bounds we are looking for more compact formulas. In any case the main meaning of the bounds is not their values, but the fact that they exist at all, hence that the last $t$ Hilbert coefficients are bounded by the first $d-t+1$ ones.

\section{Castelnuovo-Mumford regularity and Hilbert coefficients} \label{CMHil}

Let $R= \oplus_{n\ge 0}R_n$ be a Noetherian standard graded ring over a local Artinian ring $(R_0,\mm_0)$ such that $R_0/\mm_0$ is an infinite field. Let $E$ be a finitely generated graded $R$-module of dimension $d$.  For $0\le i \le d$, put
 $$ a_i(E) =
\sup \{n|\ H_{R_+}^i(E)_n \ne 0 \} ,$$
where $R_+ = \oplus _{n>  0} R_n$. The {\it Castelnuovo-Mumford regularity} of $E$  is defined by 
$$\reg(E) = \max \{ a_i(E) + i \mid  0\le  i \leq  d \},$$
and the {\it Castelnuovo-Mumford regularity of $E$  at and above level $1$}  is defined by 
$$\reg^1(E) = \max \{ a_i(E) + i \mid \  1 \leq i \leq  d \}.$$

We denote the Hilbert function $\ell_{R_0}(E_t)$ and the Hilbert polynomial of $E$ by $h_E(t)$ and $p_E(t)$, respectively. Writing $p_E(t)$ in the form:
$$p_E (t) = \sum_{i=0}^{d-1} (-1)^i e_i(E){t+d-1-i \choose d-1-i},$$
we call  the numbers $e_i(M)$  {\it the Hilbert coefficients} of $E$. 

There are different ways to bound $\reg(E)$. In this section we are interested in bounding this invariant  in terms of the Hilbert coefficients. Let $\Delta (E)$ denote the maximal generating degree of $E$. Easy examples show that one cannot bound $\reg(E)$ in terms of $\Delta (E), e_0(E),...,e_{d-1}(E)$. However these invariants bound $\reg^1(E)$, as shown in \cite[Theorem 17.2.7]{BS} and \cite[Theorem 2]{Tr1}. Below we recall the bound by Trivedi which does not depend on the number of generators of $E$ as the one in \cite{BS}. Let 
$$\Delta'(E) = \max\{\Delta(E),\ 0\}.$$
We inductively define a sequence of integers as follows: $m_1 = e_0(E) + \Delta'(E) $, and for all $i\ge 2$,
\begin{equation} \label{ENA2}
m_i = m_{i-1} + \sum_{k=0}^{i-1} (-1)^k e_k(E){m_{i-1}+i-2-k \choose i-1 -k} .
\end{equation}
Then 

\begin{Lemma} {\rm (\cite[Theorem 2]{Tr1})}  \label{TriThm} Assume that $d\ge 1$. Then $\reg^1(E) \le m_d -1$.
\end{Lemma}

The above result was originally formulated in \cite{Tr1} for $G_I(M)$, which  corresponds to the case $E$ being generated by elements of degree zero. But this assumption is not essential.  The proof  was eventually given in \cite[Lemma 4]{ST2}.  For a more algebraic proof one can use \cite[Theorem 2.7]{L1}. 

From the above bound we can derive an explicit bound for $\reg^1(E)$ in terms of $e_i(E)$ and $\Delta (E)$. However this bound is weaker.

\begin{Lemma}\label{NA2}
Let $E$ be a finitely generated  graded $R$-module of dimension $d\ge 1$.  
Put
$$\xi_{d-1}(E) = \max\{e_0(E), |e_1(E)|,...,|e_{d-1}(E)|\}.$$
 Then we have
$$\reg^1(E) \le (\xi_{d-1} (E) + \Delta'(E) + 1)^{d!} -2.$$
\end{Lemma}

\begin{pf} For short, we put $e_i := e_i(E),\ \xi := \xi_{d-1} (E)$ and $\Delta' := \Delta' (E)$. By Lemma \ref{TriThm} it suffices to show that $m_d\le (\xi + \Delta' + 1)^{d!} -1$. This is a purely arithmetic issue, which is trivial for $d=1$. By the induction hypothesis we may assume
$$m_{d-1}\le (\xi  + \Delta' + 1 )^{(d-1)!} - 1 =: \alpha  .$$
Note that
$$\sum_{i=0}^{d-1} (-1)^i e_i{\alpha +d-2-i\choose d-1-i} 
\le \xi  \sum_{i=0}^{d-1} {\alpha +d-2-i\choose d-1-i} = \xi {\alpha +d-1\choose d-1}.$$
Hence, by the recurrence formula (\ref{ENA2}) applied to $i=d$, we get
$$m_d \le \alpha + \xi {\alpha +d-1\choose d-1}.$$
If $d=2$, then $\alpha  = \xi +\Delta' $, and
$$ m_2 \le \xi + \Delta' + \xi (\xi + \Delta' + 1)  = (\xi +\Delta' + 1 ) (\xi + 1) -1\le (\xi +\Delta' + 1 )^2 - 1.$$
Assume  $d\ge 3$. Observing that ${\alpha +d-1\choose d-1} \le (\alpha+1)^{d-1}$ for all $\alpha \ge 1$ and $\alpha \ge (\xi + 1) ^2 > \xi +1$, we obtain
$$m_d \le \alpha + \xi (\alpha+1)^{d-1} \le  (1 + \xi )(\alpha+1)^{d-1} - 1 \le  (\alpha+1)^d - 1 =  (\xi +\Delta' + 1 )^{d!} -1.$$
\end{pf}

We need some more notations and definitions.  Let ($A,\mm$) be a Noetherian local  ring with an infinite residue field $K:= A/\mm$ and $M$ a finitely generated $A$-module.   Given a proper ideal $I$,  a chain of submodules
$$\FilM:\ M = M_0 \supseteq M_1 \supseteq M_2 \supseteq \cdots \supseteq M_n \supseteq  \cdots $$
 is called an {\it $I$-filtration} of $M$ if $IM_i \subseteq  M_{i+1}$ for all $i$, and  a {\it good $I$-filtration} if $IM_i =  M_{i+1}$ for all sufficiently large $i$. A module $M$ with a good $I$-filtration is called an  {\it $I$-well filtered module} (see \cite[III 2.1]{Bour}). If $N$ is a submodule of  an $I$-well filtered module $M$, then the sequence $\{M_n+N/N\}$ is a good $I$-filtration of $M/N$ and will be denoted by $\FilM/N$.

In this paper we always assume that $I$ is an $\mm$-primary ideal and $\FilM$ is a good $I$-filtration. The {\it associated graded module}  to the filtration $\FilM$ is defined by
$$G(\FilM) = \bigoplus _{n\geq 0}M_n/M_{n+1} .$$
We also say that $G(\FilM)$ is the associated module of the filtered module $M$. This is a finitely generated graded module over the standard graded ring $G:= G_I(A) := \oplus_{n\ge 0}I^n/I^{n+1}$ (see \cite[Proposition III 3.3]{Bour}). In  particular, when $\FilM$ is the $I$-adic filtration $\{I^nM\}$, $G(\FilM)$ is just the usual associated graded module $G_I(M)$. 

 We call
$H_{\FilM}(n) = \ell(M/M_{n + 1})$
the Hilbert-Samuel function of $M$ w.r.t  $\FilM$. This function agrees with  a polynomial - called the Hilbert-Samuel polynomial  and denoted by $P_{\FilM}(n)$ - for $n \gg 0$.  If we write
$$P_\FilM (t) = \sum_{i=0}^d (-1)^i e_i(\FilM){t+d-i \choose d-i},$$
then the integers $e_i(\FilM)$ are called {\it the Hilbert coefficients} of $\FilM$ (see \cite[Section 1]{RV}). When $\FilM = \{I^nM\}$, $H_{\FilM}(n)$ and $P_{\FilM}(n)$ are usually denoted by $H_{I,M}(n)$ and $P_{I,M}(n)$, respectively, and $e_i(\FilM) = e_i(I,M)$. Note that $e_i(\FilM)=e_i(G(\FilM))$ for $0\le i \le d-1$. 

Now we want to derive a bound for  $\reg(G(\FilM))$ in terms of Hilbert coefficients.   Using Lemma \ref{NA2} we  can already bound $\reg^1(G(\FilM))$ in terms of $e_0(\FilM), ..., e_{d-1}(\FilM)$. If $\depth (M)>0$, by \cite[Lemma 1.8]{DH}, $\reg(G(\FilM)) = \reg^1(G(\FilM))$, and so it is bounded in terms of $e_i(\FilM)$, $i<d$. The following example shows that this is not true if $\depth(M)=0$.

\begin{Example} \label{Ex1} {\rm Let $A= K[[x,y]]/(x^2,xy^s),\ s\ge 1$. Then $G_\mm(A) \cong k[x,y]/(x^2,xy^s)$. Since $(x^2,xy^s)$ is a so-called stable ideal, $\reg(G_\mm(A)) = s$ can be  arbitrarily large, while $e_0(A) = 1$.}
\end{Example}

Our first goal is to show that also using $e_d(\FilM)$ we can bound $\reg (G(\FilM))$. For that we need some more preparations.
We denote $M/H^0_\mm(M)$ by  $\overline{M}$ and the filtration $\FilM/ H^0_\mm(M)$ of $\overline{M}$ by $\overline{\FilM}$ and 
 let $h^0(M) = \ell(H^0_\mm(M))$.  Then

\begin{Lemma} {\rm(\cite[Proposition 2.3]{RV})} \label{NB1}
For all $n$ we have
$$ h^0(M) = P_{\FilM} (n) - P_{\overline{\FilM}}(n) = (-1)^d [e_d(\FilM) - e_d(\overline{\FilM}) ].$$
\end{Lemma}

Applying the Grothendieck-Serre formula to $G(\FilM)$  and the arguments in the proof of \cite[Lemma 3.4]{L1}, we get

\begin{Lemma}\label{NB2a}
$P_{\FilM}(n) = H_{\FilM}(n)$  for all $n\ge \reg( G(\FilM))$.
\end{Lemma}

\begin{Lemma}\label{NB2}
$h^0(M)    \le P_{\FilM}(n)$  for all $n\ge \reg( G(\overline{\FilM}))$.
\end{Lemma}

\begin{pf} By Lemma \ref{NB2a},  $P_{\overline{\FilM}}(n) = H_{\overline{\FilM}}(n)$ for all $n\ge \reg(G(\overline{\FilM}))$. Hence, by Lemma \ref{NB1},
$h^0(M)  = P_{\FilM} (n) - P_{\overline{\FilM}}(n) = P_{\FilM} (n) - H_{\overline{\FilM}}(n) \le P_{\FilM} (n) $ for all $n\ge \reg(G(\overline{\FilM}))$.
\end{pf}
 
 We call
$$r(\FilM) = \min \{r \geq 0 \mid   M_{n+1} = IM_n \ \ \text{for all} \ \ n \geq  r \}$$
the reduction number of $\FilM$ (w.r.t. $I$).   
Then we have

 \begin{Lemma}\label{A1} {\rm (\cite[Lemma 1.9]{DH})}
$\reg(G(\FilM)) \leq  \max \{ \reg(G(\overline{\FilM}));\ r(\FilM)\} + h^0(M).$
\end{Lemma}

 In the sequel we will often use the following notation:
 $$\xi_s(\FilM) = \max\{e_0(\FilM), |e_1(\FilM)|,...,|e_s(\FilM)|\},$$
 where $0\le s \le d$.  Now we can state and prove the first bound on  $\reg(G(\FilM))$ in terms of Hilbert coefficients. 
 
 \begin{Theorem} \label{NB3}  Let  $\FilM$ be  a  good $I$-filtration of $M$ of dimension $d\ge 1$.  
Then
$$\reg(G(\FilM)) < (\xi_d(\FilM)+ r(\FilM) + 1 )^{d d! + 1}  -2.$$
\end{Theorem}

\begin{pf}Let $r= r(\FilM)$, $e_i = e_i(\FilM)$ and $\xi := \xi_d(\FilM)$. 
 By \cite[Lemma 1.8]{DH}  we have 
$\reg(G(\overline{\FilM})) = \reg^1(G(\overline{\FilM})).$  By Lemma \ref{A1},
\begin{equation} \label{ENB3}
\reg(G(\FilM)) \le \max\{ \reg^1(G(\overline{\FilM})),\ r\} + h^0(M)  .
\end{equation}
Set   $\alpha := (\xi +r +1)^{d! } - 2 \ge r$.  By Lemma \ref{NB1},  $e_i(G(\overline{\FilM})) = e_i(\overline{\FilM}) = e_i$ for all $i\le d-1$. As mentioned above, $G(\overline{\FilM})$ is generated by elements of degrees at most $ r(\FilM) \ge 0$. Therefore, by Lemma \ref{NA2}, 
$\reg^1(G(\overline{\FilM})) \le \alpha $. Using (\ref{ENB3}) and  Lemma \ref{NB2} we then get
$$\begin{array}{ll}
\reg(G(\FilM)) &\le \alpha + P_\FilM(\alpha ) \\
& \le \alpha + \xi \sum_{i=0}^d{\alpha +d-i\choose d-i} \\
&= \alpha + \xi {\alpha +d+1\choose d} \\
& =  (\xi +r + 1)^{d! } - 2 + \xi {(\xi +r + 1)^{d! } -1  +d\choose d} \\
&\le (\xi +r + 1)^{d! } - 2 +  \xi  (\xi +r + 1)^{d d! } \\
& <  (\xi +r + 1)^{d d! +1 } - 2 .
\end{array}$$
\end{pf}

The above bound is a huge number when $d\gg 1$. In the case of $I$-adic filtrations of an one-dimensional  module there is a sharp bound given in a recent paper \cite{D}.

 Our next  goal is to show that  in order to bound $\reg(G(\FilM))$  one can use $\xi_{d-t}(\FilM)$,  where $t= \depth M$. For this we need some more auxiliary results.

An element $x \in I$ is called $\FilM$-{\it superficial element} for $I$ if there exists a non-negative integer $c$ such that $(M_{n+1}:_M x) \cap M_c = M_n$ for every $n \geq c$ and we say that a sequence of elements $x_1, ..., x_r$ is an $\FilM$-{\it superficial sequence} for $I$ if, for $i = 1, 2, ..., r$, $ x_i$ is an $\FilM/(x_1,...,x_{i-1})M$-superficial sequence for $I$  (see \cite[Section 1.2]{RV}). Note that $x \in I\setminus I^2$ is an $\FilM$-superficial element for $I$ if and only if its initial form $x^* \in G$ is a filter-regular element on  $G(\FilM)$, i.e. $[0 :_{G(\FilM)} x^*]_n = 0$ for all $n \gg 0$.

\begin{Lemma} \label{B1}
Let    $x$  be an $\FilM$-superficial element for $I$. Then
$$\reg(G(\FilM/xM)) \leq \reg(G(\FilM)).$$
\end{Lemma}
\begin{pf}
We have the following exact sequence:
\begin{equation} \label{EB1}
 0\longrightarrow \bigoplus_{n\ge 0}\frac{xM \cap M_n}{xM_{n-1} + xM \cap M_{n+1}} \longrightarrow  G(\FilM)/x^*G(\FilM) \longrightarrow G(\FilM/xM)\longrightarrow 0.
 \end{equation}
By \cite[Lemma 1.3(ii)]{DH} (see also \cite[Lemma 4.4]{T}), $xM \cap M_{n} = xM_{n-1}$ for $n\gg 0$. Hence
$$\reg(G(\FilM/xM) ) \leq\reg(G(\FilM)/x^*G(\FilM)) \leq \reg(G(\FilM)).$$
\end{pf}

\begin{Lemma} \label{B3}
Let $x_1, x_2, ..., x_d$ be an $\FilM$-superficial sequence for $I$. Set $M_i = M/(x_1,...,x_i)M$ and $\FilM_i = \FilM/(x_1,...,x_i)M$, where $M_0 = M$ and $\FilM_0 = \FilM$. Then, for all $0 \leq i \leq d-1$, we have
$$h^0(M_i) \leq (i+1)\xi_d(\FilM) (\reg(G(\FilM)) +2)^d \leq d \xi_d(\FilM) (\reg(G(\FilM)) +2)^d.$$
\end{Lemma}

\begin{pf}
Set $a: = \reg(G(\FilM))$ and $\xi := \xi_d(\FilM)$. We proceed by induction on $i$. Note  by Lemma \ref{B1} that $\reg(G(\overline{\FilM_i})) \leq \reg(G(\FilM_i)) \leq \reg(G(\FilM)) = a.$

For $i = 0$, by Lemma \ref{NB2}, we have
$$h^0(M_0) = h^0(M)  \leq P_{\FilM}(a)  \leq\xi \sum^d_{j = 0}{d+a-j \choose a-j} 
                = \xi {a+d+1 \choose d} \leq \xi(a+2)^d.$$
For $0 < i \leq d-1$, by \cite[Proposition 1.2]{RV}, we have $e_j(\FilM_i) = e_j(\FilM_{i-1})$ for all $0 \leq j \leq d-i-1$ and
\begin{eqnarray}
|e_{d-i}(\FilM_i)| & = & |e_{d-i}(\FilM_{i-1}) + (-1)^{d-i}\ell(0:_{M_{i-1}}x_i) | \nonumber \\
                             & \leq & |e_{d-i}(\FilM_{i-1})| + h^0(M_{i-1}) \nonumber \\
&\le & \xi_{d-i+1}(\FilM_{i-1}) + h^0(M_{i-1})  \nonumber \\
& \leq & \xi + h^0(M_{i-1}). \label{EB3}
\end{eqnarray}
Hence, by Lemma \ref{NB2}  and the induction hypothesis we get
$$\begin{array}{lll}
h^0(M_i) & \leq P_{\FilM_i}(a)  & (\text{by Lemma \ref{NB2}})\\
                & \leq \xi\sum^{d-i-1}_{j = 0}{d-i+a-j \choose d-i-j} + |e_{d-i}(\FilM_i)|  & \\
                & = \xi{a+d-i+1 \choose d-i} -\xi + |e_{d-i}(\FilM_i)|  & \\
                & \leq \xi(a+2)^{d-i} - \xi + \xi + h^0(M_{i-1})  &  (\text{by (\ref{EB3}}) \\
                & \leq \xi(a+2)^{d-i} +i\xi (a+2)^d   &  (\text{by the induction hypothesis})\\
& \leq (i+1)\xi(a+2)^d. & 
\end{array}$$
\end{pf}

\begin{Lemma} \label{B4}
Set $B = \ell(M/(x_1,x_2,...,x_d)M)$, where $x_1, x_2, ..., x_d$ is an $\FilM$-superficial sequence for $I$.  Then
$$B < (d+1)\xi_d(\FilM)(\reg(G(\FilM))+2)^d.$$
\end{Lemma}

\begin{pf}
Keep the notation in the proof of the previous lemma. Since $\dim(M_{d-1}) = 1, M_{d-1}$ is a generalized Cohen-Macaulay module. By \cite[Lemma 1.5]{CST},
$$B - e_0(x_d; M_{d-1}) = \ell(M_{d-1}/x_d M_{d-1}) - e_0(x_d; M_{d-1}) \leq h^0(M_{d-1}).$$
Since $e_0(x_d; M_{d-1}) = e_0(x_1, ..., x_d; M) = e_0(\FilM) = e_0,$ we get
 $$B \leq e_0 + h^0(M_{d-1}) \leq \xi + h^0(M_{d-1}).$$
 By Lemma \ref{B3}, $h^0(M_{d-1}) \leq d \xi(a+2)^d$ .  From this  estimate we immediately get $B < (d+1)\xi(a+2)^d$.
\end{pf}

Finally we can state and prove the second bound on  $\reg(G(\FilM))$, which only uses the first $d-t+1$ Hilbert coefficients.

\begin{Theorem} \label{B5}
Let $\FilM$ be a good $I$-filtration of $M$ with  $\dim(M) = d \geq 1$ and $\depth(M) = t$.    Then
$$\reg(G(\FilM)) \le  (\xi_{d-t}(\FilM) + r(\FilM) +1 )^{2 (d-t+1) d!} -2.$$
\end{Theorem}

\begin{pf}  For short we write $e_i = e_i(\FilM),\ \xi_s := \xi_s(\FilM)$ and $r := r(\FilM)$. We do induction on $t$. The case $t=0$ follows from Theorem \ref{NB3}. 

Assume that $t\ge 1$.
In the case $t=d$, i.e. $M$ is a Cohen-Macaulay module, the statement follows from the following bounds given in  \cite[Theorem 1.5]{DH}:
$$\reg(G(\FilM)) \le \begin{cases} e_0 +r -1 \ & \text{if}\ d=1,\\
(e_0+r+1)^{3(d-1)! -1} -d \ &  \text{if}\ d>1.
\end{cases}$$
Let $t<d$, and so $d\ge 2$. The first part of the following arguments uses the idea of the proof of  \cite[Theorem 3.3]{RTV} (see  also \cite[Theorem 4.4]{L1} and \cite[Theorem 1.5]{DH}). Let $x= x_1, ..., x_d$ be an $\FilM$-superficial sequence for $I$. Let $N= M/xM$ and $\FilN = \FilM/xM$.  Then
$\dim N = d-1$ and $\depth N = t-1$. By \cite[Proposition 1.2]{RV}, $e_i(\FilN) = e_i$ for all $i\le d-1$. Hence $\xi_{d-t}(\FilN) = \xi_{d-t}$. It is clear that $r(\FilN) \le r$. Let $m$ be an integer such that
$$m \ge \max\{r,\ \reg(G(\FilN)) \}.$$
From the exact sequence (\ref{EB1}) it follows that
$$\reg^1(G(\FilM)/x^* G(\FilM)) = \reg^1(G(\FilN)) \le m.$$
Hence, by \cite[Theorem 2.7]{L1}, 
$$\reg^1(G({\FilM})) \leq m + p_{G({\FilM})}(m).$$
By \cite[Lemma 1.6]{DH} and \cite[Lemma 1.7(i)]{DH},
$$p_{G({\FilM})}(m)  \leq H_{I,N}(m) \leq {m+d-1 \choose d-1 }\ell(N/(x_2,...,x_d)N) 
                                                =  B{m+d-1 \choose d-1}. $$
Since $\reg(G({\FilM})) = \reg^1(G({\FilM}))$ (see \cite[Lemma 1.8]{DH}),
\begin{equation} \label{EB5a}
\reg(G(\FilM))  \leq  m +  B{m+d-1 \choose d-1} \le   m+  B(m+1)^{d-1} < (B+1)(m+1)^{d-1} .
\end{equation}
Let  $\FilM_t := \FilM/(x_1,...,x_t)\FilM$. Then $\dim (M_t) = d-t$. Again, by  \cite[Proposition 1.2]{RV}, $e_i(\FilM_t) = e_i$ for all $i\le d-t$,  which yields  $\xi_{d-t}(\FilM_t) = \xi_{d-t}$.  Let $a_t : = \reg(G(\FilM_t))$.  Applying Theorem \ref{NB3} to $\FilM_t$, we have
$$a_t \le  \omega^{ (d-t) (d-t)! + 1} -2,$$
where $\omega = \xi_{d-t} + r +1$. Note that $\omega^{d-1} \ge 2^{d-1} \ge d$. Since $t\ge 1$, applying
Lemma \ref{B4} to $M_t$ we get
\begin{eqnarray}
B  & = & \ell(M_t/(x_{t+1},...,x_d)M_t)   \nonumber \\
   & \leq  &(d-t+1)\xi_{d-t} (a_t + 2)^{d-t}    \nonumber\\
& <  &d \omega  \omega^{ (d-t)^2 (d-t)! + d-t}  \nonumber\\
 & \le  & \omega^{ (d-t)^2 (d-t)! + 2d-t}.   \label{EB5b}
\end{eqnarray}
 
We distinguish two cases.

\noindent {\bf Case 1}: $t=1$. Then $\depth N = 0$. By Theorem \ref{NB3} we can take $m = \omega^{(d-1)(d-1)! +1} -2$. 
By (\ref{EB5a}) and (\ref{EB5b}) we get
$$\begin{array}{ll}
\reg(G(\FilM))  & \le  \omega^{ (d-1)^2 (d-1)! + 2d-1} ( \omega^{(d-1)(d-1)! +1} -1)^d \\
&\le  \omega^{ (d-1)^2 (d-1)! + (d-1) d!  + 3d-1} -2 \\
&\le  \omega^{ 2 d d!} - 2.
\end{array}$$

\noindent {\bf Case 2}: $t>1$. Then $d\ge 3$. By the induction hypothesis we can take $m =  \omega^{2 (d-t+1) (d-1)!} -2$. Again, by  (\ref{EB5a}) and (\ref{EB5b}) we obtain
$$\begin{array}{ll}
\reg(G(\FilM)) & < \omega^{ (d-t)^2 (d-t)! + 2d-t} (\omega ^ { 2 (d-t+1)  (d-1)!} - 1)^{d-1}\\
&\le   \omega^{ (d-t)^2 (d-t)! + 2d-2 + 2 (d-1)  (d-t+1)   (d-1)! }  - 2.
\end{array}  $$
We have 
$$2 (d-t+1)  d!=  2  (d-1) (d-t+1)  (d-1)!  +  2(d-t+1)  (d-1)!  .$$
Since $d>t \ge 2$,  the following hold
$$\begin{array}{ll}
2  (d-t+1)   (d-1)! & \ge 2 (d-1)(d-t+1)  (d-t)!\\
&> 2(d-1)  (d-t)! + 2  (d-t)^2  (d-t)! \\
&> 2(d-1) +  (d-t)^2  (d-t)!.
\end{array}$$
Hence $\reg(G(\FilM)) \le  \omega ^{2 (d-t+1) d!} -2$, as required.
\end{pf}

\begin{Remark} {\rm  Keep the notation of Lemma \ref{B3} and  Lemma \ref{B4}.  
Set 
$$B(M) = \ell (M/(x_1,...,x_d)M) \ \text{ and } \  \ 
\kappa(M) = \max\{h^0(M_i)| \ 0\leq i \leq d-1\}.$$
In the first version of this paper  (see http://viasm.edu.vn/2012/05/preprints-2012,  Preprint ViAsM12.25) we proved that 

\rm{(i)} \ $\reg(G(\FilM) )\leq B(M) + \kappa(M) + r(\FilM) - 1 \ {\mathrm{if}} \  d = 1$,

\rm{(ii)} $\reg(G(\FilM)) \leq [B(M) + \kappa(M) + r(\FilM) + 1]^{3(d-1)!-1} - d\  {\mathrm{if}}\  d\geq 2$.

Note that $B(M) = B(M_t)$ and $\kappa(M) = \kappa(M_t)$, where $t= \depth M$. Using this result,  Lemma \ref{B3}, Lemma \ref{B4}  and Theorem \ref{NB3} we can get another bound for $\reg(G(\FilM))$ in terms of $\xi_{d-t}$, which is smaller than the one of Theorem \ref{B5} if  $d-t$ is very small (compared with $d$). However, when $d-t$ is big, the bound presented in Theorem \ref{B5} is better.
} \end{Remark}

\section{Relations between Hilbert coefficients} \label{Hilb.c}

In this section we always assume that $M$ is an $A$-module of positive dimension $d$ and  $\FilM$ is a good $I$-filtration of $M$, where $I$ is an $\mm$-primary ideal. First we give an upper  bound for $(-1)^{i-1}e_i(\FilM)$ in terms of the preceding Hilbert coefficients. The first statement of the following theorem is implicitly contained in \cite{RV}.

\begin{Theorem}  \label{B0}
{\rm (i)} $e_1(\FilM) \le {e_0(\FilM) \choose 2}$.

{\rm (ii)} Let   $\xi_{i-1} := \xi_{i-1}(\FilM)$. For $i\ge 2$ we have
$$ (-1)^{i-1}e_i(\FilM) \le \xi_{i-1} {(\xi_{i-1} +r + 1)^{i!} + i \choose i} < (\xi_{i-1} +r + 1)^{i i! +1}  .$$
\end{Theorem}

\begin{pf} We do induction on $d$. Let $d=1$. Then  the inequality $e_1(\FilM) \le {e_0(\FilM) \choose 2}$ follows from \cite[Proposition 2.8 and Lemma 2.3]{RV}. 

Assume that $d\ge 2$. First we prove the statement for $i\le d-1$. Let $\overline{\FilM} = \FilM/H^0_\mm(M)$. Since $e_j(\FilM) = e_j( \overline{\FilM})$  for all $j\le d-1$, we may assume that $\depth M>0$. Let $x$ be an $\FilM$-superficial element for $I$. Then $\dim (M/xM)= d-1$ and by \cite[Proposition 1.2]{RV}, $e_j(\FilM) = e_j(\FilM/xM)$ for all $j\le d-1$. Hence, the inequalities follow from the induction hypothesis  applied to $\FilM/xM$.

Finally let $i=d$.  Since $G(\overline{\FilM})$ is generated by elements of degrees at most $r(\FilM)\ge 0$, by \cite[Lemma 1.8]{DH} and  Lemma \ref{NA2} we have
$$\reg(G(\overline{\FilM})) = \reg^1(G(\overline{\FilM})) \le (\xi_{d-1} +r + 1)^{d!} -2 =: \alpha .$$
By Lemma \ref{A1}  and Lemma \ref{NB2}  we then get
$$\begin{array}{ll}
\reg(G(\FilM)) &\le \max\{ \reg^1(G(\overline{\FilM})),\ r\} + h^0(M) \\
& \le  \max\{ \reg^1(G(\overline{\FilM})),\ r\} + P_{\FilM}(\alpha )\\
&\le  \alpha + \sum_{i=0}^{d-1} e_i(\FilM){\alpha + d-i\choose d-i} + (-1)^de_d(\FilM)\\
&\le \xi_{d-1}[\alpha -1 + \sum_{i=0}^{d} {\alpha + d-i\choose d-i}] + (-1)^de_d(\FilM) \\
& = \xi_{d-1}[\alpha -1 +  {\alpha + d +1 \choose d}] + (-1)^de_d(\FilM) \\
& < \xi_{d-1}[{\alpha + d +1 \choose d -1} +  {\alpha + d +1 \choose d}] + (-1)^de_d(\FilM) \\
&= \xi_{d-1}{\alpha + d+2\choose d} + (-1)^de_d(\FilM).
\end{array}$$
Note that ${a+d \choose d} < a^d$ for all $a\ge 4$ and $d\ge 2$. Since $\reg(G(\FilM)) \ge  r(\FilM) \ge 0$, we therefore get 
$$ \begin{array}{ll} 
(-1)^{d-1}e_d(\FilM) & \le \xi_{d-1} {\alpha + d+2\choose d} \\
& = \xi_{d-1} {(\xi_{d-1}+r + 1)^{d!}  + d\choose d} \\
& <  \xi_{d-1} (\xi_{d-1}+r +1)^{d d!} \\
& < (\xi_{d-1}+r + 1)^{d d! +1} .
\end{array}$$
\end{pf}

\begin{Remark}{\rm
Using Lemma  \ref{TriThm} and induction one can derive a better bound for  $(-1)^{i-1}e_i(\FilM)$, $i\le d-1$.  Since this bound is of almost the same complexity as the one in the above theorem, we do not give it here. The fact, that $(-1)^{i-1}e_i (I,A)$ is bounded above by a function depending on $e_0(I,A),...,e_{i-1}(I,A)$, if $i\le d-1$, was mentioned in \cite[Remark 3.10]{Bl}, provided that $A$ is an equicharacteristic local ring. Also no explicit bound was given there. }
\end{Remark}
 It is easy to see that in general $|e_i(\FilM)|$ is not bounded above by $\xi_{i-1}(\FilM)$  (see Examples \ref{Ex2} below).  In order to prove the main result of this paper, we also need bounds on  Hilbert coefficients in terms of the Castelnuovo-Mumford regularity.
 
 \vskip 0.3cm
\noindent {\it  Remark: The following result was published in the original but should be removed; see Corrigendum.}

\begin{Proposition}\label{NB5}Let $x_1,\ldots ,x_d \in I$ be an $\FilM$-superficial  sequence for $I$ and $B= \ell(M/(x_1,...,x_d)M)$. Then

{\rm (a)} For all $1\le i \le d-1$, $|e_i(\FilM)| \leq B(\reg^1(G(\FilM) )+ 1)^i$;

{\rm (b)} $|e_d(\FilM) |\le B(d+1)(\reg(G(\FilM))+1)^d$.
\end{Proposition}

\begin{pf} (a) The inequalities in (a) immediately follow from \cite[Theorem 4.6]{CHH} by noticing that $\reg(\overline{G(\FilM)}) = \reg^1(G(\FilM))$ and that $G(\FilM)$ is generated in non-negative degrees. In fact, the proof of \cite[Theorem 4.6]{CHH}  is based on \cite[Theorem 4.5(ii)]{CHH}. In its turn,  \cite[Theorem 4.5(ii)]{CHH} follows from \cite[Theorem 4.2]{CHH} and by local duality. These results were formulated for graded modules over a polynomial ring over a field. However, with a small modification,  one can show that \cite[Theorem 4.5(ii)]{CHH} and therefore also \cite[Theorem 4.6]{CHH} remain true  for any polynomial ring over an Artinian local ring.  There is yet another way:  in order to  show  \cite[Theorem 4.5(ii)]{CHH} for the case of Artinian local ring one can rewrite the proof of \cite[Theorem 4.2]{CHH} in terms of local cohomology modules. This was done in \cite[Theorem 4.1.3]{Ha}. For convenience of the reader we sketch the proof here.

{\it Claim}: Let $E$ be a graded $R$-module of dimension $d$ and let $s=\reg (E)$. Assume that $y_1,...,y_d \in R_1$ is an $E$-filter-regular sequence of $R$, that is $[0:_{E/(y_1,...,y_{i-1})E} y_i]_n =0$ for all $n\gg 0$. Put $h^i_E(t) := \ell_{R_0}(H^i_{R_+}(E)_t)$.  Our immediate aim is to show that for all  $i\ge 1$ and $s'\ge s$ we have
\begin{equation} \label{ENB5}
h^i_E(t) \le {s'-1-t \choose i-1}h_{E/(y_1,...,y_{i-1})E} (s').
\end{equation}
Since $h^i_E(t) =0$ for all $t\ge s$, we may assume  that $t\le s-1 \le s'-1$.  We proceed by induction on $i$. For $i=1$, let $E' := \oplus_{n\ge s'}E_n$. Then $\reg(E') = s'$ and  $y_1$ is regular on $E'$. The exact sequence
$$H^0_{R_+}(E'/y_1E')_t \rightarrow H^1_{R_+}(E')_{t-1} \rightarrow  H^1_{R_+}(E')_{t} \rightarrow H^1_{R_+}(E'/y_1E')_{t}$$
implies
$$h^1_{E'}(t-1) - h^1_{E'}(t)\le h^0_{E'/y_1E'} (t) \le h_{E'/y_1E'}(t).$$
Hence
$$h^1_{E}(t) \le h^1_{E'}(t) = \sum_{i=t+1}^{s'}(h^1_{E'}(i-1) - h^1_{E'}(i) ) \le \sum_{i=t+1}^{s'} h_{E'/y_1E'}(i) = h_{E'}(s') = h_E(s').$$
The case $i\ge 2$ follows from the induction hypothesis and the inequality
$$h^i_{E}(t-1) - h^i_{E}(t)\le h^{i-1}_{E/y_{i-1}E} (t).$$
So, the proof of the claim (\ref{ENB5}) in completed. Now, taking $s'=s$ and using \cite[Lemma 4.4(i)]{CHH}, we obtain
$$h^i_E(t) \le \ell (E/(y_1,...,y_d)E) {\reg(E) - 1 - t \choose i-1}{ \reg(E) + d-i \choose d-i}.$$
This is similar to the inequality  in \cite[Theorem 4.5(ii)]{CHH} and it is exactly the inequality applied in the proof of \cite[Theorem 4.6]{CHH} in order to derive (a).

(b) Let $a= \reg(G(\FilM))$ and $e_i= e_i(\FilM)$. By Lemma \ref{NB2a}, $H_\FilM(a) = P_\FilM(a)$. By \cite[Lemma 1.7]{DH}, 
$$H_\FilM(a) = \ell(M/M_{a+1}) \le \ell(M/I^{a+1}M) \le B{a+d\choose d}.$$
Since ${a+j \choose j} \le (a+1)^j$ and $\sum_{i=0}^{d}(-1)^i e_i{a+d - i\choose d-i} = H_\FilM(a)$, by (a) we get
$$\begin{array}{ll}
|e_d | &\le H_\FilM(a) + \sum_{i=0}^{d-1}|e_i|{a+d - i\choose d-i}\\
&\le B{a+d\choose d} + B\sum_{i=0}^{d-1}{a+d - i\choose d-i}(a+1)^i\\
&\le B(a+1)^d + B\sum_{i=0}^{d-1}(a+1)^{d-i}(a+1)^i \\
&= B(d+1)(a+1)^d.
\end{array}$$
\end{pf}

The following theorem is the main result of this paper:
\vskip 0.3cm
\noindent {\it Remark:  This a corrected proof. Compared with the original one, there is a small modification   in order to get the estimation (\ref{EB6}). Namely, instead of Proposition 2.3 we now use Proposition B in the Corrigendum and we correct some misprints in the estimation of  $ |e_d|$.}

\begin{Theorem} \label{B6}
Let $\FilM$ be a good $I$-filtration of $M$. Assume that $\dim(M) = d \geq 1$ and $\depth(M) = t \geq 1$. Then $|e_d(\FilM)|, |e_{d-1}(\FilM)|, ..., |e_{d-t+1}(\FilM)|$ are bounded by a function depending only on $d, e_0(\FilM), e_1(\FilM), ...,e_{d-t}(\FilM)$ and $r(\FilM)$. Namely, for all $j\ge d-t+1$ we have 
$$|e_j(\FilM)| \le  (\xi_{d-t}(\FilM) + r(\FilM) +1)^{3j (d+1-t) j!}.$$
\end{Theorem}

\begin{pf} As usual, we write $e_i := e_i(\FilM),\  \xi_s= \xi_s(\FilM)$,  $r= r(\FilM)$ and $\omega = \xi_{d-t} + r + 1$.
The case $t=d$, i.e. $M$ is a Cohen-Macaulay module, follows from the following bound given  in \cite[Theorem 1.10]{DH}:
$$|e_j| \le  (e_0 + r+1)^{3j! - j + 1} .$$
 So we can assume that $t<d$ and $d\ge 2$. 
First we prove our claim in the case  $i = d$. Let $x_1, ..., x_d$ be an $\FilM$-superficial sequence for $I$.  Keep the notation of the proof of Theorem \ref{B5}. Then by  (\ref{EB5b}) we get
$$ B \leq \omega ^{(d-t)^2  (d-t)! + 2d-t} .$$
Note that $\omega \ge 2$. Hence, by Proposition B in the Corrigendum and Theorem \ref{B5} we have
$$\begin{array}{ll}
|e_d|  &  \leq B 2^d (\reg(G(\FilM)) + 1)^d  \\
                   & < \omega ^{(d-t)^2  (d-t)! + 3d-t} \omega ^{2d (d-t+1) d!} \\
                   &  < \omega ^{(d-t)^2  (d-t)! + 2d (d-t+1) d! + 3d }.
  \end{array}$$ 
  Since $t<d$ and $d\ge 2$, it holds
  $$\begin{array}{ll}
d (d-t+1) d! & \ge 2 (d-t+1) d! \ge  2 (d-t) d (d-1)! + 2d  \\
&>  (d-t)^2  (d-t)!  + d + 2d. 
\end{array}$$
Hence
\begin{equation} \label{EB6}
|e_d| < \omega^{3 d (d-t+1) d!  }. 
\end{equation}
Now let $d-t+1 \leq j \leq d-1$. Since $\depth(M) = t$, by \cite[Proposition 1.2]{RV}, $e_j(\FilM) = e_j(\FilM_{d-j})$. Note that $\dim(M_{d-j}) = j$, $\depth(M_{d-j}) = t + j - d \geq 1$ and $r(\FilM_{d-j}) \leq r(\FilM)$. Therefore
$\xi_{d-t}(\FilM_{d-j}) = \xi_{d-t}(\FilM) = \xi_{d-t}.$
Applying (\ref{EB6}) to $\FilM_{d-j}$, we then get
$$|e_j(\FilM)| = |e_j(\FilM_{d-j})| < \omega^{3j (d-t+1)j !}.$$
\end{pf}

For the $I$-adic filtration $\FilM = \{I^nM\}_{n\geq0}$ we have $r(\FilM) = 0$. Hence, as an immediate consequence of Theorem \ref{B6}, we get the following extension of  \cite[Theorem 1]{ST2} to the non-Cohen-Macaulay case. In the Cohen-Macaulay case our bound is much bigger than that of  \cite[Theorem 1.10]{DH} (see also \cite[Theorem 4.1]{RTV} and \cite[Theorem 1]{ST2}).

\begin{Corollary}\label{B7}
Assume that $\dim(M) = d \geq 1$ and $\depth(M) = t \geq 1$. Then for all $d-t+1 \leq j \leq d$ we have
$$|e_j(I,M)| < (\xi_{d-t} + 1)^{3j (d-t+1) j!},$$
where
$$\xi_{d-t} = \max\{e_0(I,M), |e_1(I,M)|, ..., |e_{d-t}(I,M)|\}.$$
In other words, if $d-t+1 \le j \le d$, $|e_j(I,M)|$ is bounded in terms of $d, e_0(I,M),$ $ e_1(I,M), ..., e_{d-t}(I,M)$.
\end{Corollary}

Finally we can state and prove a  result about  the finiteness of Hilbert-Samuel functions.
\begin{Theorem}\label{B8}
Let $d \geq t \geq 0$, $e_0, ..., e_{d-t}$ be positive  integers. Then there exists only a finite number of Hilbert-Samuel functions associated to $d$-dimensional modules $M$ and $\mm$-primary ideals $I$ such that $\depth(M) = t$ and $e_j(I,M) \le e_j$ for all $0 \leq j \leq d-t.$
\end{Theorem}

\begin{pf}
By Corollary \ref{B7}, there exists only a finite number of Hilbert-Samuel polynomials $P_{I,M}(n)$ such that $e_j(I,M) \le e_j$ for all $0 \leq j \leq d-t.$ By Lemma \ref{NB2a}, $H_{I,M}(n) = P_{I,M}(n)$ for $n \geq \reg(G_I(M)) =: a.$ By Theorem \ref{B5},  $a$ is bounded in terms of $e_0,e_1, ..., e_{d-t}$ and $d$. Since $H_{I,M}(n) = 0$ for $n<0$ and $H_{I,M}(n)$ is an increasing function for  $n\geq0$, $H_{I,M}(n) \leq P_{I,M}(a)$ for all $n \leq a$. This implies that the number of these functions is bounded in terms of $e_0, e_1, ..., e_{d-t}$ and $d$.
\end{pf}

\begin{Example}\label{Ex2}{\rm  The following examples  show  that one cannot reduce the number of ``independent" coefficients in Theorem \ref{B6}.

\begin{itemize} 
\item[(i)] Let $A= K[[x_1,...,x_{d+1}]]/(x_1^2, x_1x_2,...,x_1x_d, x_1x_{d+1}^s)$, where $s\ge 1$, and $I = \mm = (\bar{x_1}, ..., \bar{x}_{d+1})$. Then $\dim A = d$, $\depth A = 0, \ e_0 = 1$, $e_1= \cdots = e_{d-1} = 0$, while $e_d = (-1)^d s$.

\item[(ii)] Even under certain additional assumption on $A$ we  cannot  reduce the number of ``independent" coefficients. For example, in \cite{ST1} there were constructed a complete regular local ring $R$ and an  infinite sequence  of prime ideals $\pp_n$ of $R$ such that $\dim(R/\pp_n) = 2$, $e_0(R/\pp_n) = 4$, but $e_1(R/\pp_n) = 8 -n$. 
\end{itemize}}
\end{Example}

\section*{Acknowledgment}  The authors would like to thank the referee for his/her careful reading and the list of suggested corrections which improved the presentation of this paper.

\newpage

\begin{center}
{\bf CORRIGENDUM }
\end{center}

\section*{   }

Unfortunately there was a gap in the proof of Proposition 2.3 and we have to delete it. Keeping the notation in the above  original version, then the proof  of Proposition 2.3 only gives the  following result.

\vskip0.5cm
\noindent {\bf Proposition A}.  {\it  Assume that $y_1,...,y_d \in R_1$ is an $E$-filter-regular sequence of $R$, that is $[0:_{E/(y_1,...,y_{i-1})E} y_i]_n =0$ for all $n\gg 0$. Put  $B^* =  \ell_{R_0}(E/(y_1,...,y_d)E)$. Then  
$|e_i(E)| \leq B^*(\reg^1(E)+ 1)^i$, for all $1\le i \le d-1$.}
\vskip0.5cm

These inequalities could be useful elsewhere. For the local case we can only prove 

\vskip0.5cm
\noindent {\bf Proposition B}.  {\it  Let $x_1,\ldots ,x_d \in I$ be an $\FilM$-superficial  sequence for $I$ and $B= \ell(M/(x_1,...,x_d)M)$. Then
  $|e_i(\FilM)| < B(2\reg(G(\FilM) )+ 2)^i$ for all $1\le i \le d$.}

\begin{pf}
We do induction on $d$. Let $a= \reg(G(\FilM))$ and $e_i= e_i(\FilM)$. By  Lemma 1.5, 
$$H_\FilM(a) = P_\FilM(a)=  \sum_{i=0}^{d}(-1)^i e_i{a+d - i\choose d-i}.$$
By \cite[Lemma 1.7]{DH}, 
$$H_\FilM(a) = \ell(M/M_{a+1}) \le \ell(M/I^{a+1}M) \le B{a+d\choose d}.$$
Note that ${a+j \choose j} \le (a+1)^j$ and $e_0 = e_0(I,M) \le B$. 

If $d=1$, then 
$$|e_1| = |H_\FilM(a)  - e_0(a+1)| \le \max\{B(a+1) , e_0(a+1)\} = B(a+1).$$
 
Let $d\ge 2$. First we prove the statement for $0<i\le d-1$.  Assume that $\depth(M) > 0$.  Then $\dim (M/x_1M)= d-1$ and by \cite[Proposition 1.2]{RV}, $e_i(\FilM) = e_i(\FilM/x_1M)$ for all $i\le d-1$. By Lemma 1.9, $\reg(\FilM/x_1M) \le a$. Hence, by the induction hypothesis  applied to $\FilM/x_1M$ and the sequence $x_2,...,x_d$, we get
$$|e_i(\FilM)| < B(2\reg(G(\FilM/x_1M) )+ 2)^i \le  B(2 a+ 2)^i.$$
We now assume that $\depth(M) = 0$.
Let $\overline{M} = M/H^0_\mm(M)$ and $\overline{\FilM} = \FilM/H^0_\mm(M)$. Note that $e_i\FilM) = e_i( \overline{\FilM})$  for all $i\le d-1$ and $\ell( \overline{M}/ (x_1,...,x_d)\overline{M}) \le B$. In the proof of \cite[Lemma 1.9]{DH}, it was shown that there is an exact sequence
$$0 \rightarrow  K \rightarrow  G(\FilM) \rightarrow G(\overline{\FilM}) \rightarrow 0,$$
where $K$ has a finite length. Hence $\reg (G(\overline{\FilM})) \le \reg G(M) = a$, and
$$|e_i(\FilM)| = e_i( \overline{\FilM}) < \ell( \overline{M}/ (x_1,...,x_d)\overline{M}) (2 \reg (G(\overline{\FilM})) + 2)^i \le B(2 a + 2)^i.$$
Finally, we have
$$\begin{array}{ll}
|e_d | &\le H_\FilM(a) + \sum_{i=0}^{d-1}|e_i|{a+d - i\choose d-i}\\
&<  B{a+d\choose d} + B\sum_{i=0}^{d-1} 2^i (a+1)^i  {a+d - i\choose d-i}\\
&\le B(a+1)^d + B\sum_{i=0}^{d-1} 2^i (a+1)^i  (a+1)^{d-i}\\
&= B 2^d (a+1)^d.
\end{array}$$
\end{pf}

Using Proposition B instead of Proposition 2.3 in the original  proof of  Theorem 2.4 we can still derive the same bound, because there we used a very rough estimation $d+1 \le \omega^{d+1}$, and now instead of it we only need to use the estimation $2^d \le \omega^d$. Also note that there were some misprints in establishing the inequality (8) in the original proof of  Theorem 2.4, but the inequality is correct. All these remarks were taken into account in the above corrected version.

\end{document}